\documentclass[11pt]{article}
\usepackage{amssymb,amsmath,amsthm,fullpage,mathtools,enumitem}

\usepackage{url}
\usepackage{tikz}
\usepackage{ytableau}

\usepackage{todonotes}
\newtheorem{theorem}{Theorem}

\newtheorem{lemma}{Lemma}

\theoremstyle{definition}

\newtheorem*{definition}{Definition}
\newtheorem*{notation}{Notation}

\DeclarePairedDelimiterX\makeset[2]{\{}{\}}{#1\; \delimsize\vert\; #2}
\newcommand{\pdecomp}[2]{\left[ \begin{array}{c} #1 \\ #2 \end{array}  \right]}

\newcommand{\len}{\ell}
   
\begin{document}
\title{On criteria for rook equivalence of Ferrers boards}

\author{
  Jonathan Bloom\\
  \texttt{Lafayette College}\\
  \texttt{bloomjs@lafayette.edu}\\
  \and 
  Dan Saracino\\
    \texttt{Colgate University}\\
    \texttt{dsaracino@colgate.edu}
}

\maketitle
\begin{abstract} In \cite{BloomSaracino:Rook2017} we introduced a new notion of Wilf equivalence of integer partitions and proved that rook equivalence implies Wilf equivalence.  In the present paper we prove the converse and thereby establish a new criterion for rook equivalence. We also refine two of the standard criteria for rook equivalence and establish another new one involving what we call \emph{nested sequences of L's}.
\end{abstract}

\maketitle

\section{Introduction}
In \cite{BloomSaracino:Rook2017} we introduced a notion of Wilf equivalence of integer partitions (Ferrers boards) and proved that rook equivalence implies Wilf equivalence.  The main purpose of the present paper (see Section ~\ref{sec:wilf}) is to prove the converse and thereby establish a new criterion for rook equivalence that is quite different from the known criteria.

Of course our result also provides a characterization of Wilf equivalence, and we point out that characterizations of Wilf equivalence in other contexts have remained conjectures despite serious attempts to find proofs.  One notable example is the so-called ``Rearrangement Conjecture", which involves a notion of Wilf equivalence on words (see \cite{KitaevLieseRemmelSagan:Rational2009,PantoneVatter:On-the-R}), and another is Conjecture 1 in~\cite{AlbertBouvel:A-genera2014}, involving Catalan structures.  

During our study of Wilf equivalence of integer partitions, we also discovered refinements of two standard criteria for rook equivalence.   Foata and Sch\"{u}tzenberger proved in~\cite{foata1970rook} that two Ferrers boards are rook equivalent if and only if one can be obtained from the other via a series of transformations called $(i,j)$-transformations.  (We recall the definition in Section~\ref{sec:transforms}.)  We prove, in Theorem~\ref{thm:rook iff transform}, that $(i,1)$-transformations suffice.  Foata and Sch\"{u}tzenberger also established the well-known criterion (recalled in Section~\ref{sec:Ls}) for rook equivalence in terms of multisets of integers obtained from the rows of a Ferrers board.  We refine this criterion in Section~\ref{sec:Ls} by showing (Theorem~\ref{thm:rook iff same salient multiset2}) that the multisets obtained from certain special rows (which we call \emph{salient}) suffice to determine rook equivalence.  Finally, using this refinement, we give (in Theorem~\ref{thm:nested2}) another new criterion for rook equivalence in terms of what we call \emph{nested sequences of $L$'s}.    

Throughout this paper we regard a partition $\mu$ of $n$ as an infinite weakly decreasing sequence of nonnegative integers $\mu_1\geq \mu_2\geq \cdots$ whose nonzero terms sum to $n$.    We call $n$ the \emph{weight} of $\mu$ and write $|\mu|=n$. We refer to $\mu_i$ as the \emph{$i$th part} of $\mu$.  We define the \emph{length} $\len(\mu)$ of $\mu$ to be the number of nonzero parts.  In the special case that all the nonzero parts are distinct, we call $\mu$ \emph{strict}.

Via the following theorem of Foata and Sch\"{u}tzenberger, strict partitions play a recurring role in what follows.
\begin{theorem}(Foata, Sch\"{u}tzenberger \cite[Theorem 11]{foata1970rook})\label{thm:FS strict}
Every partition $\mu$ is rook equivalent to a unique strict partition.  
\end{theorem}

\section{Rook equivalence and Wilf equivalence}\label{sec:wilf}

We begin by recalling several definitions from~\cite{BloomSaracino:Rook2017}.

Viewing integer partitions as Ferrers boards, we say that a partition $\alpha$ \emph{contains (as a pattern)} a partition $\mu$ if $\mu$ can be obtained by deleting rows and columns of $\alpha$. For example, if $\mu=(4,3,3, 2,2)$ then by deleting the rows and columns indicated in red we see that $\mu$ is contained in the partition
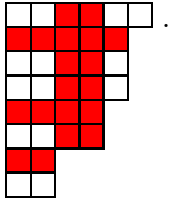
\ytableausetup{boxsize=.8em}
$$\ydiagram{6,5,5,5,4,4,2,2}*[*(red)]{0,5,0,0,4,0,2}*[*(red)]{2+2,2+2,2+2,2+2,2+2,2+2}\ .$$

Denoting by $\mathcal{P}_n(\mu)$ the set of all partitions of $n$ that contain $\mu$,  we say that partitions $\mu$ and $\nu$ are \emph{Wilf equivalent} if, for every integer $n\geq 1$, $|\mathcal{P}_n(\mu)|=|\mathcal{P}_n(\nu)|$.  In the argument that follows we will need to refine $\mathcal{P}_n(\mu)$ by defining, for every integer $k\geq 0,$ 
$$\mathcal{P}_n(\mu,k)=\makeset{\alpha\in \mathcal{P}_n(\mu)}{ \alpha_1=k+\mu_1}.$$

\begin{notation}

 For partitions $\mu$ and $\nu,$ we write $\mu\sim_r \nu$ to indicate that $\mu$ and $\nu$ are rook equivalent, and we write $\mu\sim_W \nu$ to indicate that $\mu$ and $\nu$ are Wilf equivalent.
\end{notation}

We will rely on the following theorem, which was the main theorem of~\cite {BloomSaracino:Rook2017}.

\begin{theorem}(Bloom, Saracino \cite[Theorem 3]{BloomSaracino:Rook2017})\label{thm:rook implies wilf}
If $\mu$ and $\nu$ are partitions of $n$ and $\mu\sim_r \nu$, then $\mu\sim_W \nu.$
\end{theorem}

Our first objective is to prove the following theorem.  

\begin{theorem}\label{thm:rook iff Wilf} If $\mu$ and $\nu$ are partitions of $n$, then $\mu\sim_r \nu$  if and only if $\mu\sim_W \nu.$.
\end{theorem}

By Theorem~\ref{thm:rook implies wilf}, we need only prove that Wilf equivalent partitions are rook equivalent. To prove this, it suffices to show that two Wilf equivalent strict partitions must be identical. For suppose $\pi$ and $\rho$ are any two partitions such that $\pi\sim_W \rho.$  By  Theorem~\ref{thm:FS strict} we have  $\pi\sim_r \pi'$ and $\rho\sim_r\rho',$ where $\pi'$ and $\rho'$ are strict. By Theorem~\ref{thm:rook implies wilf} we then have $\pi\sim_W \pi'$ and $\rho\sim_W \rho'$. Therefore $\pi'\sim_W \rho'$, so, assuming that this forces $\pi'=\rho',$ we have $\pi\sim_r \pi'\sim_r \rho,$  so $\pi\sim_r \rho$.

Now let $\mu$ and $\nu$ be two distinct strict partitions of $n$, with the aim of proving that $\mu\nsim_W \nu.$ Let $r$ be the largest integer such that $\mu_r\neq \nu_r$. Note that $r> 1$ since $|\mu|= |\nu|$. Without loss of generality, we assume that $\mu_r< \nu_r\neq 0$, and we let $m=r+\nu_r$.  We prove that $\mu\nsim_W \nu$ by showing that $|\mathcal{P}_{n+m-1}(\mu)|< |\mathcal{P}_{n+m-1}(\nu)|$. To establish this inequality, we show that 

\begin{equation}\label{eq:converse goal}
  |\mathcal{P}_{n+m-1}(\mu,k)|\leq |\mathcal{P}_{n+m-1}(\nu,k)|  
\end{equation}
for all $k\geq 0$, and that the inequality is strict when $k=1$.

Our first lemma will facilitate the comparison of $\mathcal{P}_{n+m-1}(\mu,k)$ and $\mathcal{P}_{n+m-1}(\nu,k)$.

\begin{lemma}\label{lem:top rows} We have 
$\nu_1 \geq m-1$ and $\mu_1 > m-1$.
\end{lemma}

\begin{proof}
Since $\nu$ is strict and we have chosen $\nu_r\neq 0$, we see that $1+\nu_1\geq r+\nu_r=m$, so $\nu_1 \geq m-1$.

To see that $\mu_1 > m-1$, note that since $|\mu|=|\nu|$ and $\mu_i=\nu_i$ for all $i> r$, we have

\begin{equation}\label{eqn:sums} 
\sum_{i=1}^r (i+\mu_i)= \sum_{i=1}^r (i+\nu_i).
\end{equation}
Since $\nu$ is strict and $\nu_r\neq 0$, the right side of (\ref{eqn:sums}) is at least $r(r+\nu_r)=rm$. Now consider the summation on the left and write it as 
\begin{equation}\label{eqn:sums2}
    \sum_{i=1}^r (i+\mu_i) = \sum_{i=1}^s (i+u_i) + \sum_{i=s+1}^r i,
\end{equation}
where $s$ is such that $\mu_i>0$ for $1\leq i\leq s$ and $\mu_i = 0$ for $s+1\leq i\leq r$.  
For a contradiction, assume that $\mu_1\leq m-1$.  As $\mu$ is strict it follows that $i+\mu_i \leq m$ for $1\leq i\leq s$.     Now consider the terms in the rightmost sum of (\ref{eqn:sums2}).  By definition of $m$ we see that 
$r\leq r+\mu_r< r+\nu_r=m$ and therefore each term (if any) in this sum must be less than $m$.  So each term in the left side of (\ref{eqn:sums}) is at most $m,$ and $r+\mu_r< m.$ Therefore the left side of (\ref{eqn:sums}) is less than $rm$, contradicting our previous observation that the right side of (\ref{eqn:sums}) is at least $rm$.  We conclude that $\mu_1>m-1$.
\end{proof}

To further our comparison of the sets $\mathcal{P}_{n+m-1}(\mu,k)$ and $\mathcal{P}_{n+m-1}(\nu,k)$ we next consider building the elements in these two sets by inserting rows and columns into the partitions $\mu$ and $\nu$, respectively.  In particular, we make use of  Lemma 1 of~\cite{BloomSaracino:Rook2017}, which states that any element $\gamma\in \mathcal{P}_{n+m-1}(\nu,k)$ can be obtained by first adding to $\nu$ (as new columns) the columns of a partition $\alpha$ with $\alpha_1 = k$ (obtaining a partition which we denote by $\nu+\alpha$), and then adding (as new rows) the rows of a partition $\beta$ with $\beta_1\leq \nu_1+k$, such that $|\alpha|+|\beta|=m-1$. We refer to $\gamma$ as the \emph{extension} of $\nu$ by the pair $(\alpha,\beta)$.  

In our situation any pair $(\alpha,\beta)$  with $\alpha_1 = k$ and $|\alpha|+|\beta|  = m-1$ can be used to construct an element in the aforementioned sets since Lemma~\ref{lem:top rows} and the fact that
$$\beta_1 \leq |\alpha|+|\beta| = m-1$$
guarantee that $\beta_1\leq \mu_1, \nu_1$.  

In light of this, we would now like to show that for any two pairs of partitions $(\alpha, \beta)$ and $(\alpha', \beta')$ such that $\alpha_1=\alpha'_1 = k$ and each pair has cumulative weight $m-1$, the following must hold: if the extensions of $\nu$ obtained from these pairs coincide, then the extensions of $\mu$ obtained from these pairs must also coincide.  We show this in Lemma~\ref{lem:pairs coincide}, but first we need a preliminary lemma.

% By Lemma~\ref{lem:top rows} above, every $\beta$ such that $|\alpha|+|\beta|=m-1$ has $\beta_1\leq \nu_1$, so we can use all $\beta$'s such that $|\alpha|+|\beta|=m-1$. A similar argument applies to elements of $\mathcal{P}_{n+m-1}(\mu,k)$, so to  prove that $|\mathcal{P}_{n+m-1}(\mu,k)|\leq |\mathcal{P}_{n+m-1}(\nu,k)|$ it suffices to show that for any two pairs of partitions $(\alpha, \beta)$ and $(\alpha', \beta')$ such that $\alpha_1=\alpha'_1 = k$ and each pair has cumulative weight $m-1$, the following must hold: if the extensions of $\nu$ obtained from these pairs coincide then the extensions of $\mu$ obtained from these pairs must also coincide.

\begin{lemma}\label{lem:unique construction} Let  $\nu$ be a strict partition and fix $1\leq a \leq \len(\nu)$. Consider two pairs of partitions $(\alpha,\beta)$ and $(\alpha',\beta')$ with $\len(\alpha) < a$ and $\beta_1 < \nu_a$.  If the extensions of $\nu$ by these two pairs coincide, then $(\alpha,\beta)=(\alpha',\beta')$.
\end{lemma}
% \begin{proof}
% 	Assume $\gamma$ is the extension of $\nu$ by the pair $(\alpha,\beta)$.  As  $\len(\alpha)<a$,  $\beta_1<\nu_a$ and $\nu$ is strict, it follows that $\gamma$ has exactly $a$ parts of size at least $\nu_a$ and, furthermore, that $a$ is the only row of $\gamma$ such that $\gamma_a = \nu_a$.  Now if $\gamma$ is also obtained by extending $\nu$ by the pair $(\alpha',\beta')$ it follows that $\len(\alpha') < a$ and $\beta'_1 < \nu_a$.  From this it is easy to see that $\beta' = \beta$ and, in turn, that $\alpha' = \alpha$.  
% \end{proof}

\begin{proof}
	Assume $\gamma$ is the extension of $\nu$ by the pair $(\alpha,\beta)$.  As  $\len(\alpha)<a$,  $\beta_1<\nu_a$ and $\nu$ is strict, it follows that $\gamma_a=\nu_a$ and that $\nu$ and $\gamma$ each have exactly $a$ parts of size greater than or equal to $\nu_a.$
	
	Now assume $\gamma$ is also obtained by extending $\nu$ by the pair $(\alpha',\beta')$. Then $\len(\alpha')< a$ because $\gamma_a=\nu_a$, and $\beta_1< \nu_a$ because $\nu$ and $\gamma$ each have exactly $a$ parts of size greater than or equal to $\nu_a.$ This implies that $(\alpha,\beta)= (\alpha',\beta'),$ as we see by considering separately the rows in $\gamma$ above/below row $a.$
\end{proof}

For the remainder of this section let $\mu$, $\nu$, $r$, and $m$ be as defined in the paragraph preceding Lemma~\ref{lem:top rows}.

\begin{lemma}\label{lem:pairs coincide}
	  Fix pairs of partitions $( \alpha, \beta)$ and $(\alpha',\beta')$ such that $\alpha_1 = \alpha'_1=k$ and both pairs have total weight $\nu_r + (r-1)$.  If the extensions of $\nu$ by  these pairs coincide, then the extensions of $\mu$ by these pairs must also coincide.  

In the case $k=1$, if $\len(\alpha) > \len(\alpha')$ and both pairs yield the same extension of $\nu$, then $ \len(\alpha),\len(\alpha')\geq r$.

\end{lemma}

\begin{proof}
 To start, we make a general observation.  For any positive integers $x$ and $y$ with 
$$x + y = \nu_r + (r-1)\quad \text{and}\quad x\leq r-1$$
we claim that $y \leq \nu_{x+1}$. In fact, since $\nu$ is strict and $x\leq r-1$ we have
$$\nu_{x+1} - \nu_r \geq r-(x+1).$$
This immediately implies that $\nu_{x+1} \geq \nu_r + (r-1) -x = y$ as claimed.  

The lemma is obvious when $k=0$ so we assume $k\geq 1$.  We now prove the first claim of the lemma by considering two cases.  

	\medskip
	Case 1: $|\alpha| \leq r-1$. 
	\medskip
	
	First consider the subcase when $k>1$.  In this subcase, it is obvious that $\len(\alpha) <|\alpha|$. By the above general observation applied to $x = |\alpha|$ and $y = |\beta|$ we also see that 
$$\beta_1 \leq |\beta| \leq \nu_{|\alpha|+1}< \nu_{|\alpha|}.$$
  It now follows from Lemma~\ref{lem:unique construction} (with $a= |\alpha|$) that $(\alpha,\beta) = (\alpha',\beta')$, which yields the first claim of our lemma.  

When $k=1$ it will suffice to show that $|\alpha| = |\alpha'|$, i.e., $\alpha = \alpha'$, from which it follows (since the extensions of $\nu$ by $(\alpha,\beta)$ and $(\alpha',\beta')$ coincide) that $\beta = \beta'$, proving the first claim of our lemma in this case.  We can assume, for a contradiction, that $|\alpha| <|\alpha'|$, for if $|\alpha'|< |\alpha|$ then $|\alpha'|\leq r-1$   and we can interchange the roles of $\alpha$ and $\alpha'$.  From our general observation, with $x= |\alpha|$, we see that 
\begin{equation}\label{eq:width beta}
\beta_1 \leq |\beta| \leq \nu_{|\alpha| + 1} = \nu_{x + 1}.
\end{equation}
Denoting the multiset of parts of  $\nu + \alpha$ as
$$\{\nu_1+1 > \nu_2 +1 > \cdots > \nu_x +1> \nu_{x+1} > \cdots\},$$
it follows that since $\nu$ is strict this multiset has no part of size $\nu_{x+1}+1$.  On the other hand, recalling our assumption that $|\alpha|<|\alpha'|$, consider the largest $x+1$ elements in the multiset of parts for $\nu+\alpha'$:
$$\{\nu_1+1 > \nu_2 +1 > \cdots > \nu_x +1> \nu_{x+1}+1 \}.$$
Observe that it does contain a part of size $\nu_{x+1}+1$.  Consequently, $\beta$ must have a part of size $\nu_{x+1}+1$ as we demand that the extension of $\nu$ by $(\alpha,\beta)$ coincide with the extension by $(\alpha',\beta')$. Since (\ref{eq:width beta}) forbids this as an option, we reach our contradiction.  We conclude that $|\alpha| = |\alpha'|$ as needed.  
	
\medskip
Case 2:  $|\alpha| > r-1$.
\medskip

By the proof of Case 1 we can assume that $|\alpha'| > r-1$, as otherwise we obtain the contradiction that $\alpha' = \alpha$ with $|\alpha'| \leq r-1$.  Since the total weight of each pair is $\nu_r + (r-1)$  we can conclude that  $|\beta| < \nu_r$ and $|\beta'| < \nu_r$.  Therefore the largest $r$ parts of the extension of $\nu$ by $(\alpha,\beta)$ are 
$$\nu_1+ \alpha_1 > \nu_2+ \alpha_2 > \cdots > \nu_r+ \alpha_r.$$
Since the extension of $\nu$ by $(\alpha,\beta)$  coincides with the extension of $\nu$ by  $(\alpha',\beta')$, these $r$ largest parts are also
$$\nu_1+ \alpha_1' > \nu_2+ \alpha_2' > \cdots > \nu_r+ \alpha_r'.$$
Hence we may conclude that $\alpha_i = \alpha_i'$ for $i \leq r$.  Recalling that our definition of $r$ insists that $(\nu_{r+1}, \nu_{r+2} , \ldots ) = (\mu_{r+1}, \mu_{r+2}, \ldots )$, we conclude that the first claim of the lemma holds in case 2 as well. 

The previous two cases establish the first claim of our lemma.  To address the second claim, observe that if $(\alpha,\beta)$ and $(\alpha',\beta')$ yield the same extension of $\nu$ and  $k=1$, then the proof of Case 1 guarantees that $\len(\alpha) = \len(\alpha')$ whenever  $\len(\alpha)=|\alpha| \leq r-1$ or $\len(\alpha')=|\alpha'| \leq r-1$.  Hence the only way we can have $\len(\alpha) > \len(\alpha')$ is if  $ \len(\alpha), \len(\alpha')\geq r$.  This establishes our second claim.  
\end{proof}

\begin{lemma}\label{lem:k=1} 
We have $|\mathcal{P}_{n+m-1}(\mu,1)|< |\mathcal{P}_{n+m-1}(\nu,1)|$.
\end{lemma}

\begin{proof}
The first claim of the preceding lemma establishes that $|\mathcal{P}_{n+m-1}(\mu,1)|\leq  |\mathcal{P}_{n+m-1}(\nu,1)|$.  To demonstrate that this inequality is strict,  we show that there exist pairs  $(\alpha,\beta)$ and $(\alpha',\beta')$ with total weight $\nu_r + r-1$ that yield the same extension of $\mu$ but different extensions of $\nu$. (To aid the reader an example of this construction is given following the proof.)

To construct such pairs, recall that by definition of $r$ we have 
$$r+\mu_r \leq \nu_r+r -1= m-1.$$
Hence there exists some $j\geq 0$ such that 
$$r + \mu_r+j  = \nu_r+r -1= m-1.$$
Let $\alpha$ be a column of length $r$ and let $\beta$ be the partition with parts $\mu_r$ and $j$ so that the total weight of the pair $(\alpha,\beta)$ is $m-1$ as needed.  Likewise, define $\alpha'$ to be a column of length $r-1$ (recall that $r> 1$) and let $\beta'$ be the partition with parts $\mu_r+1$ and $j$ so that $(\alpha',\beta')$ also has total weight $m-1$.  Then $(\alpha,\beta)$ and $(\alpha',\beta')$ yield the same extension of $\mu$ because adding column $\alpha$ and a row of size $\mu_r$ to $\mu$ yields the same result as adding column $\alpha'$ and a row of size $\mu_r+1$. On the other hand, the second claim in the preceding lemma guarantees that $(\alpha,\beta)$ and $(\alpha', \beta')$ do not yield the same extension of $\nu$, since $\len(\alpha')<r$.
\end{proof}

To illustrate the coinciding extensions of $\mu$  described in the proof of Lemma~\ref{lem:k=1} consider the case when $\nu = (7,6,5)$ and $\mu = (10,7,1)$ so that $r = 3$ and $j = 3$.   We now have 

$$ \ydiagram{11}*[*(red)]{1+1}*
{0,8}*[*(red)]{0,1+1}*
[*(blue)]{0,0,3}*
{0,0,0,1}*[*(red)]{0,0,0,1+1}*
[*(blue)]{0,0,0,0,1}
\qquad\text{ and }\qquad 
\ydiagram{11}*[*(red)]{1+1}*
{0,8}*[*(red)]{0,1+1}*
[*(blue)]{0,0,3}*
{0,0,0,0,1}*
[*(blue)]{0,0,0,2}$$
where the partition on the left is the extension of $\mu$ by  $\alpha = (1,1,1)$, shown in red, and $\beta = (3, 1)$, shown in blue.  The partition on the right is the extension of $\mu$ by the pair $\alpha' = (1,1)$ and $\beta' = (3, 2)$. Note that the extensions of $\nu$ by $(\alpha, \beta)$ and $(\alpha', \beta')$ are different because they have smallest parts $1$ and $2$, respectively.

\begin{proof}[Proof of Theorem~\ref{thm:rook iff Wilf}]
As described in the paragraph following the statement of  Theorem~\ref{thm:rook iff Wilf}, we need only show that distinct strict partitions $\mu$ and $\nu$ cannot be Wilf equivalent. Lemma~\ref{lem:pairs coincide} tells us that 
$$|\mathcal{P}_{n+m-1}(\mu,k)|\leq |\mathcal{P}_{n+m-1}(\mu,k)|$$
for all $k\geq 1$, and Lemma~\ref{lem:k=1} proves that this inequality is strict when $k=1$.  So $|\mathcal{P}_{n+m-1}(\mu)| \neq |\mathcal{P}_{n+m-1}(\nu)|$ and hence $\mu\nsim_W \nu.$   This completes our proof.
\end{proof}

\section{$(i,1)$-Transforms}\label{sec:transforms}

    Let $\mu$ be a partition (viewed as a Ferrers board), and let $(i,j)$ denote the box in the $i$th row from the top and the $j$th column from the left. If $(i,j)$ has the property that replacing the subboard $\makeset{(x,y)}{x\geq i, \ y\geq j}$ by its conjugate yields a partition, then following Foata and Sch\"{u}tzenberger in \cite{foata1970rook}, we call this new partition the \emph{(i,j)-transform} of $\mu$ and denote it by $\mu(i,j)$. For example, we have
    \ytableausetup{boxsize=.8em}
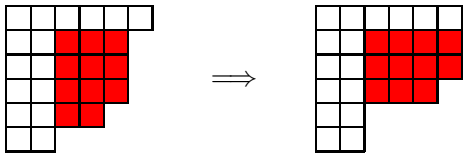
$$\ydiagram{6,5,5,5,2,2}*[*(red)]{0,2+3,2+3,2+3,2+2} \qquad\raisebox{-.75cm}{$\Longrightarrow$} \qquad \ydiagram{6,6,6,5,2,2}*[*(red)]{0,2+4,2+4,2+3}$$
where the partition on the right is the $(2,3)$-transform of the partition on the left. 

Foata and Sch\"{u}tzenberger gave the following  characterization of rook equivalence in terms of these transforms.  
\begin{theorem}(Foata, Sch\"{u}tzenberger \cite[Corollary 6, Lemma 9, Theorem 11]{foata1970rook})\label{thm:FS transforms}
    Partitions $\mu$ and $\tau$ of $n$ are rook equivalent if and only if one can be obtained from the other by a sequence of $(i,j)$-transformations.
\end{theorem}

Our purpose in this section is to establish the nonobvious fact that if partitions $\mu$ and $\tau$ are rook equivalent then one can be obtained from the other by a sequence of $(i,1)$-transformations, i.e., $\mu$ and $\tau$ are \emph{$(i,1)$-equivalent}.  

\begin{lemma}
	Every partition of $n\geq 1$ is $(i,1)$-equivalent to a strict partition.  
\end{lemma}
\begin{proof}
	We proceed by induction on $n$, the result being obvious for $n=1$.  Consider any partition $\pi$ of $n$, and assume the result for all partitions of $m$, for all $m< n$.

Choose $\sigma$ in the $(i,1)$-equivalence class of $\pi$ with the top row of $\sigma$ as large as possible. Write $\sigma= \pdecomp{\rho}{\tau}$, where $\rho$ is the top row of $\sigma$ and $\tau$ is the portion of $\sigma$ below $\rho$.  The left column of $\tau$ must be shorter than $\rho$, for otherwise the conjugate $\sigma^*$ of $\sigma$ would contradict our choice of $\sigma$.  Likewise, the top row of $\tau$ must be shorter than $\rho,$ for otherwise  $$\pdecomp{\rho}{\tau^*}^*$$ would contradict our choice of $\sigma$.

By the induction hypothesis, $\tau$ is $(i,1)$-equivalent to a strict partition, so we have a sequence of partitions $\tau,\tau^{(1)}, \ldots, \tau^{(k)}$ such that $\tau^{(k)}$ is strict and each $\tau^{(j)}$ is an $(i,1)$-transform of the partition that precedes it in the list.

Since the left column of $\tau$ is shorter than $\rho$, we see that $\pdecomp{\rho}{\tau^{(1)}}$ is obtained from  $\pdecomp{\rho}{\tau}$ by a single $(i,1)$-transformation.  By the same reasoning as above, the left column and top row of $\tau^{(1)}$ are both shorter than $\rho$.  Now, iterating this argument, we conclude that $\sigma$ (and therefore $\pi$) is $(i,1)$-equivalent to the strict partition $\pdecomp{\rho}{\tau^{(k)}}$.
\end{proof}

\begin{theorem}\label{thm:rook iff transform}
	If $\mu$ and $\nu$ are partitions of $n$, then $\mu\sim_r \nu$ if and only if $\mu$ and $\nu$ are $(i,1)$-equivalent.
\end{theorem}
\begin{proof}
	By Theorem~\ref{thm:FS transforms}, it suffices to show that if $\mu\sim_r\nu$, then $\mu$ and $\nu$ are $(i,1)$-equivalent. For any partitions $\mu$ and $\nu,$ the preceding lemma tells us that $\mu$ and $\nu$ are $(i,1)$-equivalent to strict partitions $\mu'$ and $\nu'$, respectively.  By Theorem~\ref{thm:FS transforms} and the fact that $\mu\sim_r\nu$ we have  
	$\mu' \sim_r\mu \sim_r \nu\sim_r \nu'$.
	By Theorem~\ref{thm:FS strict} we have $\mu' = \nu'$.  Therefore $\mu$ and $\nu$ are $(i,1)$-equivalent.
\end{proof}

\section{Nested Sequences of $L$'s}\label{sec:Ls}

In this final section we characterize rook equivalence in terms of what we call \emph{nested sequences of $L$'s}.

\begin{definition}
If $\mu$ is a partition, then an $L$ in $\mu$ is a set of boxes in $\mu$ obtained by starting with a box $(i,j)$, which we call the \emph{corner} of $L$, and taking all boxes $(i,x)$ with $x\geq j$ and all boxes $(y,j)$ with $y\leq i$.  A \emph{nested sequence of L's} in $\mu$ is a sequence  $L_1,L_2,\ldots, L_k$ of $L$'s in $\mu$ such that for $1\leq c<d\leq k$ the corner of $L_d$ is strictly above and strictly to the right of the corner of $L_c$.  

By the \emph{weight} of a nested sequence of $L$'s we mean the total number of boxes in all the $L$'s in the sequence.  For any partition $\mu$ and positive integer $k$, we let $w_k(\mu)$ denote the maximum weight of a sequence of $k$ nested $L$'s in $\mu$.  	
\end{definition}

For an example of a nested sequence of $L$'s consider the sequence $L_1,L_2,L_3$ as depicted below:
    \ytableausetup{boxsize=.8em}
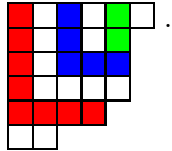
$$\ydiagram{6,5,5,5,2,2}*[*(red)]{1,1,1,1,4}*[*(blue)]{2+1,2+1,2+3}*[*(green)]{4+1,4+1}\ .$$
The weight of this sequence is 15. 

The $L'$s in our definition are reminiscent of the hooks of a Ferrers board, the obvious difference being that in a hook the vertical leg extends downward from the corner, instead of upward.  We are not aware of any known connection between rook equivalence and hooks.

One of our main results in this section is the following theorem.

\begin{theorem}\label{thm:nested2}
	Partitions $\mu$ and $\nu$ are rook equivalent if and only if $w_k(\mu) = w_k(\nu)$ for all positive integers $k$.  
\end{theorem}

To prove this result we will need a refinement of the multiset criterion for rook equivalence established by Foata and Sch\"utzenberger. To state their criterion, we associate to each partition $\mu$ a multiset $\mathcal{M}_\mu$ as follows.

\begin{definition} For each partition $\mu,$ let $\mathcal{M}_\mu$ be the multiset $\mathcal{M}_\mu=\makeset{i+\mu_i}{i\geq 1}$.
\end{definition}

\begin{theorem}(Foata, Sch\"{u}tzenberger \cite{foata1970rook})\label{thm: FS multiset}
Partitions $\mu$ and $\nu$ of $n$ are rook equivalent if and only if $\mathcal{M}_\mu=\mathcal{M}_\nu.$
\end{theorem}

The idea behind our refinement is that only those elements of the multiset corresponding to what we call the \emph{salient} rows of a partition are needed to determine rook equivalence.

\begin{definition}
	We call a row $i$ in a partition $\mu$ \emph{salient} if there exists some $j>i$ such that $i+ \mu_i\geq j+\mu_j$.
\end{definition}

For example, consider the partition $\mu = (4,3,3,2,2,2)$.  The values of $i+ \mu_i$ are 
$$(\underline{5},5,\underline{6},6,\underline{7},\underline{8},\textbf{7,8,9, \ldots})$$
where the values in bold correspond to parts of size 0 and the underlined values correspond to the salient rows.  We point out that one can easily determine the nonsalient rows as they correspond to the (strict) ``right-to-left minima" in the displayed sequence.   It is also not hard to show (but we will not need this result) that a row $s$ of $\mu$ is salient if and only if the southwest walk along the diagonal with slope 1 that starts at the rightmost box in row $s$ encounters the bottom box of some column.

\begin{definition}
    For each partition $\mu,$ let $\mathcal{S}_\mu$ and $\mathcal{N}_{\mu}$ be the multisets $$\mathcal{S}_\mu=\makeset{i+\mu_i}{ \textrm{row} \ i\  \textrm{is salient in}\  \mu}\ \textrm{and}\ \mathcal{N}_{\mu}=\makeset{i+\mu_i}{\textrm{row}\ i\ \textrm{is not salient in}\ \mu}.$$
\end{definition}

Our definition of the nonsalient rows immediately implies that $\mathcal{N}_{\mu}$ is actually a set, i.e., all its elements have multiplicity one, and that $[\len(\mu)+1,\infty)\subseteq \mathcal{N}_{\mu}.$

Our objective is to prove that in Theorem~\ref{thm: FS multiset} we can replace $\mathcal{M}_\mu$ by $\mathcal{S}_\mu.$ To do so, we will need
an additional definition, reminiscent of the well-studied Durfee square rank.  

\begin{definition}
	We define the \emph{staircase rank} of a partition $\mu$ to be the largest integer $k$ such that
	$$\mu_i\geq k-i+1\text{ for } 1\leq i\leq k.$$
\end{definition}
It is easy to see that if $n$ is the Durfee square rank of a partition and $k$ is its staircase rank, then $n\leq k$.  These ranks can be distinct, as one can see by considering the partition $(2,1)$. 

A straightforward check (whose details we leave to the reader) shows that staircase rank is invariant under $(i,j)$-transformations. So, by Theorem~\ref{thm:FS transforms}, any two rook equivalent partitions must have the same staircase rank. We make repeated use of this observation below.

\begin{lemma}\label{lem:not salient2} For a partition $\mu$ of $n$ with staircase rank $k$, we have 
$\mathcal{N}_{\mu} = [k+1,\infty)$. 
\end{lemma} 

\begin{proof} By the definition of staircase rank we have $i+\mu_i\geq k+1$ for all $i$, so $\mathcal{N}_{\mu}\subseteq [k+1,\infty)$.  As noted above, $[\len(\mu)+1, \infty)\subseteq \mathcal{N}_{\mu}.$  Now consider any $t\in [k+1, \len(\mu)]$ so that the northeast walk with slope 1 that starts at $(t-1,1)$ must contain a first box of the form $(j,\mu_j)$. Then the boxes
$$(t-1,1), (t-2,2),\cdots, (j,\mu_j)$$
are in $\mu$ as are the boxes 
$$(t-1,2), (t-2,3),\cdots, (j+1,\mu_j).$$
From this it follows that for any $j<i<t$ we have $j+\mu_j = t < i+\mu_i$.  Since $t\leq \len(\mu)$ we further see that when $i\geq t$ we have $j+\mu_j = t <i+ \mu_i$.  So row $j$ is not salient, and therefore $t\in \mathcal{N}_{\mu}.$  Since $t$ is an arbitrary element of $[k+1,\len(\mu)],$ we conclude that $[k+1,\len(\mu)]\subseteq \mathcal{N}_{\mu}$. This completes the proof.
\end{proof}

From this lemma it immediately follows that
$\mathcal{S}_{\mu}=\mathcal{M}_{\mu}\setminus [k+1,\infty).$
Since both $\mathcal{M}_\mu$ and the staircase rank are preserved under $(i,j)$-transformations the following lemma is clear.

\begin{lemma}\label{lem:salient preserved under transform2}
Let $\mu$ and $\nu$ be partitions of $n$ that differ by a sequence of $(i,j)$-transformations.  Then  $\mathcal{S}_{\mu}=\mathcal{S}_{\nu}$.

\end{lemma}

For the next two lemmas we make use of the observation that if $\mu'$ is the unique strict partition rook equivalent to $\mu$, then we have 
\begin{align}\label{eq:salient strict}
\mathcal{S}_{\mu} = \mathcal{S}_{\mu'} = \makeset{s+\mu'_s}{1\leq s\leq \len(\mu')},   
\end{align}
where the first equality follows by the previous lemma and the second equality follows since every nonzero row of a strict partition is salient.  We also observe that for every strict partition $\delta$ its length is the number of elements in $\mathcal{S}_{\delta}$ counted with multiplicity, and $\ell(\delta)$ is also the staircase rank of $\delta$.

\begin{theorem}\label{thm:rook iff same salient multiset2}
	Partitions $\mu$ and $\nu$ of $n$ are rook equivalent if and only if $ 
\mathcal{S}_{\mu} = \mathcal{S}_{\nu}.$ 
\end{theorem}

\begin{proof}
	Assuming $\mu\sim_r\nu$, Theorem~\ref{thm:FS transforms} tells us that $\nu$ can be obtained from $\mu$ through a sequence of $(i,j)$-transformations.  Lemma~\ref{lem:salient preserved under transform2} then guarantees $\mathcal{S}_{\mu}=\mathcal{S}_{\nu}$.  
	
	For the converse, assume that 
	$\mathcal{S}_{\mu}=\mathcal{S}_{\nu}.$
	Let $\mu'$ and $\nu'$  be the unique strict partitions rook equivalent to $\mu$ and $\nu$ respectively.  By (\ref{eq:salient strict}) we see that 
	$\mathcal{S}_{\mu'} =\mathcal{S}_{\mu} =\mathcal{S}_{\nu} = \mathcal{S}_{\nu'}$,
	so (by the last sentence preceding the statement of this theorem) $\len(\mu')=\len(\nu')$ and therefore $\mu'$ and $\nu'$ have the same staircase rank.  This rank is also the staircase rank of $\mu$ and $\nu$, so, since $\mathcal{S}_{\mu}=\mathcal{S}_{\nu}$, we have $\mathcal{M}_{\mu}=\mathcal{M}_{\nu}$ by the paragraph preceding Lemma~\ref{lem:salient preserved under transform2}.  Thus $\mu\sim_r \nu$ by Theorem~\ref{thm: FS multiset}.  
\end{proof}

\begin{lemma}\label{lem:staircase2}
The staircase rank of $\mu$ is the number of salient rows in $\mu$.  Furthermore, the length of the $k$th salient row from the bottom is at least $k$.  
\end{lemma}

\begin{proof}
Let $\mu'$ be the unique strict partition rook equivalent to $\mu$. 
Since the staircase rank and the number of salient rows are both invariant under rook equivalence, and $\ell(\mu')$ is both the staircase rank of $\mu'$ and the number of salient rows in $\mu'$, our first claim is proved. 

To address the second claim, suppose the $k$th salient row from the bottom is the $i$th row of $\mu$.  By deleting the rows above $i$ we may consider a partition $\tau=(\mu_i,\mu_{i+1},\ldots)$ with exactly $k$ salient rows, where the first row is salient.  By our first claim $\tau$ has staircase rank $k$ and thus $\tau_1 \geq k$.  This proves our second claim.  
\end{proof}

\begin{lemma}\label{lem:max weight in salient2}
If a sequence of $k$ nested $L$'s  in $\mu$ has weight $w_k(\mu)$, then the corner of each $L$ in the sequence resides in a salient row.  
\end{lemma} 

\begin{proof}
Let $\mathcal{L}$ be a sequence of nested $L$'s with weight $w_k(\mu)$ and denote by $r_1 < \cdots < r_k$ the rows containing the corners of these $L$'s.  Since $\mathcal{L}$ has weight $w_k(\mu)$, the corners of the $L$'s in $\mathcal{L}$ must reside in the first $k$ columns.  It will suffice to show that if the corner of one of these $L$'s were in a non-salient row then we could find a sequence $\mathcal{L}^*$ of $k$ nested $L$'s in $\mu$ with larger weight.

First consider the case when row $r_1$ is not salient and all subsequent $L$'s have their corners in salient rows.  Restricting our attention to the truncated partition $\tau = (\mu_{r_1},\mu_{r_1+1},\ldots)$ we see that  $\tau_1\geq k$ as this row intersects the vertical legs of $k$ nested $L$'s.  Likewise, $\tau_2\geq k-1$, and so on.  We conclude that $\tau$ has staircase rank at least $k$, and therefore, by Lemma~\ref{lem:staircase2}, that $\tau$ contains at least $k$ salient rows.  Hence $\mu$ contains at least $k$ salient rows below row $r_1$.  Therefore we can choose some salient row $s>r_1$ which is distinct from the rows $r_2,\ldots, r_k$.  By the definition of salience we know that $s+ \mu_s > r_1+\mu_{r_1}$.  Among these $k$ salient rows the $i$th one from the bottom has length at least $i$ by Lemma~\ref{lem:staircase2}.  Therefore there exists a sequence $\mathcal{L}^*$ of $k$ nested $L$'s having their horizontal legs on rows $s,r_2,\ldots, r_k$ and their vertical legs in the first $k$ columns.  Denoting the weights of $\mathcal{L}$ and $\mathcal{L}^*$ by $|\mathcal{L}|$ and $|\mathcal{L}^*|$,  we see that
    $$|\mathcal{L}|  = \sum_{i =1}^k (r_i+\mu_{r_i}-i)=\sum_{i=1}^k(r_i+ \mu_{r_i})-\frac{k(k+1)}{2} $$
and that $|\mathcal{L}^*|$ is obtained from this last sum by replacing the term $r_1+ \mu_{r_1}$ by $s+ \mu_s$.  Therefore $|\mathcal{L}^*|> |\mathcal{L}|$.

This completes our proof in the case where row $r_1$ is not salient and rows $r_2,\ldots, r_k$ are salient.
In all other cases we can consider the first non-salient row $r_i$ from the bottom and use the above argument to replace only the $L$'s with corners on rows $r_i$ and below.
\end{proof}

With these lemmas in hand, we now turn to the proof of Theorem~\ref{thm:nested2}.

\begin{proof}[Proof of Theorem~\ref{thm:nested2}]
	We let $w'_k(\mu)$ denote the maximum weight of a sequence of $k$ nested $L$'s in $\mu$ whose horizontal legs lie in salient rows of $\mu$.  By Lemma~\ref{lem:max weight in salient2} we know that $w_k(\mu) = w'_k(\mu)$ for all $k>0$ and $\mu$, so we aim to prove that $\mu\sim_r\nu$ if and only if $w'_k(\mu) = w'_k(\nu)$ for all positive $k$.

	Let $w''_k(\mu)$ denote the maximum weight of $k$ $L$'s in $\mu$ whose horizontal legs lie in distinct salient rows and whose vertical legs are in the first  column of $\mu$.  It is clear that $w_k'(\mu)$ is the weight of a sequence of $k$ nested $L$'s with their horizontal legs in salient rows of $\mu$ and their vertical legs in the first $k$ columns of $\mu$.  By sliding the vertical legs of these $L$'s into the first column of $\mu$ we see that $w'_k(\mu)+ \frac{(k-1)k}{2} \leq w''_k(\mu)$.  (The vertical legs slide distances $0,1,2,\ldots, k-1,$ respectively.) On the other hand, for any set of $L$'s in $\mu$ (of the type indicated in the definition of $w_k''(\mu)$) that yields the value $w''_k(\mu)$ we can, by Lemma~\ref{lem:staircase2}, move the vertical legs of these $L$'s to the right and obtain a sequence of $k$ nested $L$'s with weight $w''_k(\mu)-\frac{(k-1)k}{2}$.  So $w'_k(\mu)\geq w''_k(\mu)-\frac{(k-1)k}{2}$.  We have shown that $w'_k(\mu)+\frac{(k-1)k}{2} = w''_k(\mu)$ for all $k$ and $\mu$, so our goal now is to show that $\mu\sim_r\nu$ if and only if $w''_k(\mu) = w''_k(\nu)$ for all positive $k$.
	
	But $w''_k(\mu) = w''_k(\nu)$ for all positive $k$ if and only if we have the multiset equality
	$\mathcal{S}_{\mu}=\mathcal{S}_{\nu}$.
	By Theorem~\ref{thm:rook iff same salient multiset2}, this multiset equality holds if and only if $\mu\sim_r \nu$.  
\end{proof}

\bibliography{mybib}
\bibliographystyle{siam}      
\end{document}